\documentclass[12pt]{article}
\usepackage{amsmath}
\usepackage{amssymb,amsmath}
\def\fc{{\textbf{\textit c}}}
\def\D{\Delta}
\def\a{\alpha}
\def\b{\beta}
\def\SM{\!\setminus\!}

\def\UU{{\mathcal U}}
\def\LL{{\mathcal L}}
\def\DD{{\cal D}}
\def\Der{{\rm Der}}
\def\Inn{{\rm Inn}}
\def\Ker{{\rm Ker}}
\def\Im{{\rm Im}}

\def\cl{\centerline}
\def\rar{\rightarrow}

\def\vs{\vspace*}
\def\WW{\mathcal {W}}
\def\ni{\noindent}
\def\VV{\mathcal {V}}
\def\Z{\mathbb{Z}}

\def\F{\mathbb{F}}
\def\QED{\hfill$\Box$}

\numberwithin{equation}{section}
\newtheorem{theo}{Theorem}[section]
\newtheorem{defi}[theo]{Definition}

\newtheorem{lemm}[theo]{Lemma}
\newtheorem{prop}[theo]{Proposition}

\newtheorem{case}{Case}

\begin{document}

\cl{{\bf Lie bialgebra structures on the Schr\"{o}dinger-Virasoro
Lie algebra}\footnote {\!\!\!Supported by NSF grants 10825101,
10671027 of China, China Postdoctoral Science Foundation Grant
20080440720, Science Foundation of University Doctoral Program CNCE,
Ministry of Education of China}} \vs{6pt}

\cl{ Jianzhi Han$^{\dag)}$, Junbo Li$^{\dag,\ddag)}$, Yucai
Su$^{\dag)}$}

\cl{\small $^{\dag)}$Department of Mathematics, University of
Science\! and \!Technology \!of\! China, Hefei 230026, China}

\cl{\small $^{\ddag)}$Department of Mathematics, Changshu Institute
of Technology, Changshu 215500, China}

\cl{\small E-mail: jzzhan@mail.ustc.edu.cn,\ sd\_junbo@163.com,\
ycsu@ustc.edu.cn}\vs{6pt}

{\small\parskip .005 truein \baselineskip 3pt \lineskip 3pt

\noindent{{\bf Abstract.} In this paper we investigate  Lie
bialgebra structures on the Schr\"{o}dinger-Virasoro algebra $\LL$.
Surprisingly, we find out an interesting fact that not all Lie
bialgebra structures on the Schr\"{o}dinger-Virasoro algebra are
triangular coboundary, which is different from the related known
results of some Lie algebras related to the Virasoro
algebra.\vs{5pt}

}

}
\parskip .001 truein\baselineskip 6pt \lineskip 6pt

\vs{18pt}

\cl{\bf\S1. \
Introduction}\setcounter{section}{1}\setcounter{equation}{0}

To search for the solutions of the Yang-Baxter quantum equation,
Drinfeld \cite{D1} %originally????
introduced the notion of Lie bialgebras in 1983. Since then, many
papers on Lie bialgebras appeared, e.g., [\ref{EK}, \ref{LSX},
\ref{M1}--\ref{NT}, \ref{SS}, \ref{T}, \ref{WSS0}, \ref{WSS}].
%\cite{EK,LSX}, \cite{M1}--\cite{NT}, \cite{SS,T,WSS0,WSS}.
%\cite{LSX}--\cite{T}, \cite{WSS0,WSS}.
Witt type Lie bialgebras introduced in \cite{T} were classified in
\cite{NT}, whose generalized cases were considered in
\cite{SS,WSS0}. Lie bialgebra structures on generalized
Virasoro-like and Block  Lie algebras were %respectively
investigated in \cite{WSS, %} and \cite{
LSX}. The
Schr\"{o}dinger-Virasoro Lie algebra  \cite{H1} was introduced in
the context of non-equilibrium statistical physics during the
process of investigating the free Schr\"{o}dinger equations. There
are two sectors of this type Lie algebras, i.e., the original one
and the twisted one, both of which are closely related to the
Schr\"{o}dinger algebra and the Virasoro algebra, which play
important roles in many areas of mathematics and physics (e.g.,
statistical physics) and have been investigated in a series
of papers [\ref{GJP}, \ref{H1}, \ref{H3}--%, \ref{LS1},
\ref{LS2},
\ref{LSZ}, \ref{RU}, \ref{TZ}].
% \cite{LS1}--\cite{U}, \cite{TZ,ZT}.
However, Lie bialgebra structures on the Schr\"{o}dinger-Virasoro
Lie algebra have not yet been considered. Drinfel'd \cite{D2} posed
the problem whether or not there exists a general way to quantilize
all Lie bialgebras. Although Etingof and Kazhdan \cite{EK} gave a
positive answer to the question, they did not provide a uniform
method to realize quantilizations of all Lie bialgebras. As a matter
of fact, investigating Lie bialgebras and quantilizations is a
complicated problem. In this paper we shall determine  Lie bialgebra
structures on the Schr\"{o}dinger-Virasoro algebra $\LL$. It is
known that every Lie bialgebra structure on the Lie algebras
considered in \cite{{LSX},NT,SS,WSS} is triangle coboundary.
Surprisingly, we find out an interesting fact that not all Lie
bialgebra structures on the Schr\"{o}dinger-Virasoro algebra are
triangular coboundary.

The Schr\"{o}dinger-Virasoro algebra $\LL$  \cite{H1} is an
infinite-dimensional Lie algebra over a field $\F$ of characteristic
0 with basis $\{L_n,Y_p,M_n\,|\,n\in \Z,p\in \frac{1}{2}+\Z\}$ and
the following non-vanishing Lie brackets
\begin{eqnarray}\label{LB}\begin{array}{lll}
&&[L_m,L_{n}]=(n-m)L_{n+m},\ \ \
[L_m,M_n]=nM_{n+m},\\[6pt]
&&[\,L_n,Y_p\,]=(p-\frac{n}{2})Y_{p+n},\ \ \ \ \,\
[\,\,Y_p,Y_{q}\,\,]=(q-p)M_{p+q}.\end{array}
\end{eqnarray}
It has an infinite-dimensional ideal $S$ with basis
$\{Y_{n+\frac{1}{2}},M_n\,|\,n\in \Z\}$ and a Witt subalgebra (the
centerless Virasoro algebra) $\WW$ with  basis $\{L_n\,|\,n\in
\Z\}$. And $\F M_0$ is the center of $\LL$.

Let us recall the definitions related to Lie bialgebras. Let $L$ be
any vector space. Denote $\xi$ the {\it cyclic map} of $L\otimes
L\otimes L$, namely, $ \xi (x_{1} \otimes x_{2} \otimes x_{3})
=x_{2} \otimes x_{3} \otimes x_{1}$ for $x_1,x_2,x_3\in L,$ and
$\tau$ the {\it twist map} of $L\otimes L$, i.e., $\tau(x\otimes y)=
y \otimes x$ for $x,y\in L$. The definitions of a Lie algebra and
Lie coalgebra can be reformulated as follows. A {\it Lie algebra} is
a pair $(L,\delta)$ of a vector space $L$ and a bilinear map $\delta
:L\otimes L\rar L$ with the conditions:
\begin{eqnarray*}
&&\Ker(1-\tau) \subset \Ker\,\delta,\ \ \ \delta \cdot (1 \otimes
\delta ) \cdot (1 + \xi +\xi^{2}) =0 :  L \otimes L\otimes L\rar
 L.
\end{eqnarray*}
Dually, a {\it Lie coalgebra} is a pair $(L,\D)$ of a vector space $
L$ and a linear map $\D: L\to L\otimes L$ satisfying:
\begin{eqnarray}\label{cLie-s-s}
&&\Im\,\D \subset \Im(1- \tau),\ \ \ (1 + \xi +\xi^{2}) \cdot (1
\otimes \D) \cdot \D =0: L\to L\otimes L\otimes L.
\end{eqnarray}
We shall use the symbol ``$\cdot$'' to stand for the {\it diagonal
adjoint action}:
\begin{eqnarray*}
&&x\cdot(\mbox{$\sum\limits_{i}$}{a_{i}\otimes b_{i}})=
\mbox{$\sum\limits_{i}$}({[x,a_{i}]\otimes
b_{i}+a_{i}\otimes[x,b_{i}]}).
\end{eqnarray*}
%\begin{defi}\rm
A {\it Lie bialgebra} is a triple $( L,\delta,\D)$ such that $( L,
\delta)$ is a Lie algebra, $( L,\D)$ is a Lie coalgebra, and the
following compatible condition holds:
\begin{eqnarray}
\label{tr}&&\D\delta (x\otimes y) = x \cdot \D y - y \cdot \D x,\ \
\forall\,\,x,y\in L.
\end{eqnarray}
%\end{defi}
Denote  $\UU$ the universal enveloping algebra of $ L$, and  $1$ the
identity element of $\UU$. For any $r =\sum_{i} {a_{i} \otimes
b_{i}}\in L\otimes L$, define $\fc(r)$ to be elements of $\UU
\otimes \UU \otimes \UU$ by
\begin{eqnarray*}
&&\fc(r) = [r^{12} , r^{13}] +[r^{12} , r^{23}] +[r^{13} , r^{23}],
\end{eqnarray*}
where $r^{12}=\sum_{i}{a_{i} \otimes b_{i} \otimes 1} , \ \ r^{13}=
\sum_{i}{a_{i} \otimes 1 \otimes b_{i}} , \ \ r^{23}=\sum_{i}{1
\otimes a_{i} \otimes b_{i}}$. Obviously
\begin{eqnarray*}
\fc(r)=\mbox{$\sum\limits_{i,j}$}[a_i,a_j]\otimes b_i\otimes
b_j+\mbox{$\sum\limits_{i,j}$}a_i\otimes [b_i,a_j]\otimes b_j+
\mbox{$\sum\limits_{i,j}$}a_i\otimes a_j\otimes [b_i,b_j].
\end{eqnarray*}

\begin{defi}\label{def2}
\rm (1) A {\it coboundary Lie bialgebra} is a $4$-tuple $( L,
\delta, \D,r),$ where $( L,\delta,\D)$ is a Lie bialgebra and $r \in
\Im(1-\tau) \subset L\otimes L$ such that $\D=\D_r$ is a {\it
coboundary of $r$}, where $\D_r$ is defined by
\begin{eqnarray}
\label{D-r}\D_r(x)=x\cdot r\mbox{\ \ for\ \ }x\in L.
\end{eqnarray}

(2) A coboundary Lie bialgebra $( L,\delta,\D,r)$ is called {\it
triangular} if it satisfies the following {\it classical Yang-Baxter
Equation} (CYBE):
\begin{eqnarray}
\label{CYBE} \fc(r)=0.
\end{eqnarray}

(3) An element $r\in\Im(1-\tau)\subset L\otimes L$ is said to
satisfy the \textit{modified Yang-Baxter equation} (MYBE) if
\begin{eqnarray}
\label{MYBE}x\cdot \fc(r)=0,\ \,\forall\,\,x\in L.
\end{eqnarray}
\end{defi}

 Denote  $\VV= \LL\otimes \LL$. Then $ \LL$ and $\VV$ are both $
\frac12{\mathbb{Z}}$-graded. Denote (see \S3) $\Der(\LL,\VV)$ (resp.
$\Inn(\LL,\VV)$) the space of derivations (resp. inner derivations)
from $\LL$ to $\VV$, and $H^1(\LL,\VV)$ the first cohomology group
of $\LL$ with coefficients in $\VV$.  For any 6 elements
$\a,a^\dag,\b,\b^\dag,\gamma,\gamma^\dag\in\F$, one can easily
verify that the linear map $D:\LL\to\VV$ defined below is a
derivation\vs{-7pt}:
\begin{eqnarray}\label{def-D}&&D(L_n)=(n\a+\gamma)M_0\otimes
M_n+(n\a^\dag+\gamma^\dag))M_n\otimes M_0,\nonumber\\
&&D(Y_{n-\frac12})=\b M_0\otimes Y_{n-\frac12}+\b^ \dag
Y_{n-\frac12}\otimes M_0,\nonumber\\&&
 D(M_n)=2(\b M_0\otimes
M_n+\b^\dag M_n\otimes M_0),\ \ n\in \Z.\end{eqnarray}
 Denote $\DD$
the 6-dimensional space spanned by the such elements $D$.
 Let $\DD_1$ be the subspace of $\DD$ consisting of
elements $D$ such that $D(\LL)\subseteq \mathrm{Im}(1-\tau)$.
Namely, $\DD_1$ is the 3-dimensional subspace of $\DD$ consisting of
elements $D$ with $\a=-\a^\dag,\,\b=-\b^\dag,\,\gamma=-\gamma^\dag$.

 The main results of
this paper can be formulated as follows.
\begin{theo}\label{theo}\vskip-3pt
 \begin{itemize}\parskip-2pt
\item[\rm(i)] $\Der(\LL,\VV)=\mathrm{Inn}(\LL,\VV)\oplus\DD$ and
$H^1(\LL,\VV)=\Der(\LL,\VV)/\Inn(\LL,\VV)\cong \DD.$

\item[\rm(ii)]  Let $(\LL,[\cdot,\cdot],\D)$ be a Lie bialgebra such that $\D$
has the decomposition $\D_r+D$ with respect to
$\Der(\LL,\VV)=\mathrm{Inn}(\LL,\VV)\oplus\DD$, where
$r\in\VV\,({\rm mod\,}M_0\otimes M_0)$ and $D\in \DD$. Then, $r\in
\mathrm{Im}(1-\tau)$ and $D\in \DD_1$. Furthermore,
$(\LL,[\cdot,\cdot],D)$ is a Lie bialgebra provided $D\in\DD_1$.

\item[\rm(iii)] A Lie bialgebra $(\LL,[\cdot,\cdot],\D)$  is triangular coboundary
if and only if $\D$ is an inner derivation $($thus $\D=\D_r$, where
$r\in\mathrm{Im}(1-\tau)$ is some solution of CYBE$)$.
\end{itemize}\end{theo}

Note that Theorem \ref{theo}(ii)  shows that there exist Lie
bialgebras which are not triangular coboundary. Moreover, Theorem
\ref{theo}(iii) (resp.~(ii)) gives a description of Lie bigalebra
structures determined by $\mathrm{Inn}(\LL,\VV)$ (resp.~$\DD$). But
one cannot expect that $(\LL,[\cdot,\cdot],\D_r+D)$ would
automatically become a Lie bialgebra even if both
$(\LL,[\cdot,\cdot],\D_r)$ and $(\LL,[\cdot,\cdot],D)$ are Lie
bialgebras, since the equation in (\ref{cLie-s-s}) does not
satisfies the linear relation for derivations.  This is  one of the
reasons why it is  difficult to classify all Lie bialgebra
structures (in case when all Lie bialgebra structures are triangular
coboundary, the classification of Lie bialgebra structures is
equivalent to solving all solutions of CYBE (cf.~Lemma \ref{some}),
which is not done even for the Virasoro algebra
\cite{NT}).\vskip18pt

\cl{\bf\S2. \ Some preliminary
results}\setcounter{section}{2}\vspace*{-1pt}
\setcounter{theo}{0}\setcounter{equation}{0} Throughout the paper,
we denote by $\Z_+$ the set of all nonnegative integers and $\Z^*$
(resp. $\F^*$) the set of all nonzero elements of $\Z$ (resp. $\F$).

\begin{lemm}\label{Legr}\rm
Regard $\LL^{\otimes n}$ $($the tensor product of $n$ copies of
$\LL)$ as an $\LL$-module under the adjoint diagonal action of
$\LL$. Suppose $r\in\LL^{\otimes n}$ satisfying $x\cdot r=0$,
$\forall$ $x\in\LL$. Then $r\in\F M_0^{\otimes n}$.
\end{lemm}\vs{-6pt}
\noindent{\it Proof}\ \ It can be proved directly by using the
similar arguments as those presented in the proof of Lemma 2.2 of
\cite{WSS}.\QED
\begin{lemm}\rm \label{some}
\begin{itemize}\parskip-3pt\item[\rm(i)] $r$ satisfies CYBE in $(\ref{CYBE})$ if and
only if it satisfies MYBE in $(\ref{MYBE})$.
\item[\rm(ii)]
 Let $\LL$ be a Lie
algebra and $r\in\Im(1-\tau)\subset \LL\otimes\LL,$\  then
\begin{eqnarray}
\label{add-c}(1+\xi+\xi^{2})\cdot(1\otimes\D_r)\cdot\D_r(x)=x\cdot
\fc(r),\ \ \forall\,\,x\in\LL,
\end{eqnarray}
and the triple $(\LL,[\cdot,\cdot], \D_r)$ is a Lie bialgebra if and
only if $r$ satisfies CYBE $(\ref{CYBE})$. \end{itemize}
\end{lemm}
\noindent{\it Proof } (i)\ \ If $r$ satisfies MYBE, by Lemma
\ref{Legr}, $\fc(r)\in\F M_0^{\otimes3}$. As in \cite{NT}, one has
$\fc(r)=0$. The reverse statement is obvious.

(ii) The result can be found in \cite{D1,D2,NT}.\QED\vskip5pt

The following technical result gives some descriptions of solutions
of CYBE.
\begin{prop}\label{props}
Let $ r=\sum_{q\in\frac12\Z}r_q\in\mbox{\rm Im}(1-\tau)$ be a
nonzero solution of CYBE, and $p$ be the maximal index with $r_p\neq
0$.
 If %$\fc(r)=0$ and
$p\in\Z$, then $r_p\in\cup_{i=0}^5V_i$, where $V_1,...,V_5$ are
subspaces spanned respectively by
$$\begin{array}{llll}
L_0\otimes L_p-L_p\otimes L_0,& M_0\otimes L_p-L_p\otimes
M_0;\\[4pt]
L_0\otimes L_p-L_p\otimes L_0,& L_0\otimes M_p-M_p\otimes L_0;\\[4pt]
M_0\otimes L_p-L_p\otimes
M_0,&M_0\otimes M_p-M_p\otimes M_0;\\[4pt]
L_0\otimes M_p-M_p\otimes L_0,&M_0\otimes M_p-M_p\otimes M_0;\\[4pt]
M_j\otimes M_{p-j}-M_{p-j}\otimes M_j\mbox{ \ for all \ }j\in\Z.
\!\!\!\!\!\!\!\!\!\!\!\!\!\!\!\!\!\!\!\!\!\!\!\!\end{array}$$ If
%$\fc(r)=0$ and
$p\in\frac{1}{2}+\Z$, then $r_p\in V_6\cup
V_7\cup(\cup_{i\in\Z}V_8^{(i)})$, where $V_6,\,V_7,\,V_8^{(i)}$ are
subspaces spanned respectively by
$$\begin{array}{lll}
L_0\otimes Y_p-Y_p\otimes L_0,& M_0\otimes Y_p-Y_p\otimes
M_0;\\[4pt]
L_{\frac{2}{3}p}\otimes Y_{\frac{1}{3}p}-Y_{\frac{1}{3}p}\otimes
L_{\frac{2}{3}p},& M_{\frac{2}{3}p}\otimes
Y_{\frac{1}{3}p}-Y_{\frac{1}{3}p}\otimes M_{\frac{2}{3}p};
\\[4pt]
 M_i\otimes Y_{p-i}-Y_{p-i}\otimes M_i.\!\!\!\!\!\!\!\!\!\!\!\!\!\!\!\end{array}$$
Here and below, we treat $L_a,Y_b$ as zero if
$a\notin\Z,\,b\notin\frac12+Z$.
\end{prop}
\noindent{\it Proof } First assume $p\neq0$. We can suppose  $p>0$
otherwise the arguments are similar. Let $a_i\otimes
b_{p-i}-b_{p-i}\otimes a_i$ be a term in $r_p$ with nonzero
coefficient, where $a_i,b_i\in\{X_i,M_i,Y_{i-\frac12}\}$. Now we
prove it case by case\vspace*{-7pt}.

\begin{case}\label{caseLLYY}
$a_i\otimes b_{p-i}=L_i\otimes L_{p-i}$ or $Y_{i-\frac12}\otimes
Y_{p+\frac12-i}$ \vspace*{-9pt}.
\end{case}

\rm Changing the sign of the coefficient of $a_i\otimes
b_{p-i}-b_{p-i}\otimes a_i$ if necessary, we may assume that $i>0$,
since $p>0$. Moreover, we can assume that $i$ is maximal. Assume
that $a_i\otimes b_{p-i}=L_i\otimes L_{p-i}$. Then $r_p$ cannot
contain  terms of the form $d(X_j\otimes W_{p-j}-X_{p-j}\otimes
W_j)$ with $d\in\F^*$ and $j\neq i$, where $X,W\in\{L,M\}$. Suppose
the contrary. We could take a nonzero term $d_0(X_{j_0}\otimes
W_{p-j_0}-W_{p-j_0}\otimes L_{j_0})$ of $r_p$ with $j_0\neq i$ being
maximal, and one could easily see that $[L_i,X_{j_0}]\otimes
L_{p-i}\otimes W_{p-j_0}$ would be a term in $\fc(r)_{2p}$ with
nonzero coefficient, contradicting $\fc(r)_{2p}=0$. In particular,
we have shown that the $L_i\otimes L_{p-i}-L_{p-i}\otimes L_i$ is
the unique term of the form $L_k\otimes L_{p-k}-L_{p-k}\otimes
L_k(k\in\Z)$ in $r_p$. Thus, $0=\fc(L_i\otimes
L_{p-i}-L_{p-i}\otimes L_i)\in\LL\otimes\LL\otimes \LL$,  since the
terms of the type $L_{k_1}\otimes L_{k_2}\otimes L_{k_3}$ of
$\fc(r)_{2p}$ can only be obtained from $L_i\otimes
L_{p-i}-L_{p-i}\otimes L_i$. It follows that $i=p$, i.e.,
$d(L_p\otimes L_0-L_0\otimes L_p)$ is the only term in $r_p$ of the
form $L_k\otimes L_{p-k}-L_{p-k}\otimes L_k(k\in\Z)$. Applying the
similar arguments as above to  the case $a_i\otimes
b_{p-i}=Y_{i-\frac12}\otimes Y_{p+\frac12-i}$ one can see that
$Y_{i-\frac12}\otimes Y_{p+\frac12-i}- Y_{p+\frac12-i}\otimes
Y_{i-\frac 12}$ is the unique term of the form $Y_{k-\frac12}\otimes
Y_{p+\frac12-k}- Y_{p+\frac12-k}\otimes Y_{k-\frac 12}(k\in\Z)$.
Clearly, one also should have $\fc(Y_{i-\frac12}\otimes
Y_{p+\frac12-i}- Y_{p+\frac12-i}\otimes Y_{i-\frac 12})=0$, which is
impossible. Thus, the  situation $a_i\otimes
b_{p-i}=Y_{i-\frac12}\otimes Y_{p+\frac12-i}$ cannot occur. Whence
we conclude that $r_p$ must lie in the subspace spanned by
$X_0\otimes W_p-W_p\otimes X_0$, where $X,W\in\{L,M\}$. Furthermore,
observe that the coefficient of $M_0\otimes M_p-M_p\otimes M_0$ must
be  zero.  Thus,
$$r_p\in\,\rm{ Span}\{L_0\otimes L_p-L_p\otimes L_0, L_0\otimes
M_p-M_p\otimes L_0, M_0\otimes L_p-L_p\otimes M_0\}.$$ Now one can
check that either $r_p\in V_1$ or $r_p\in V_2$.  Namely, $r_p\in
V_1\cup V_2$. \vspace*{-7pt}.

\begin{case}\label{caseLM}
$a_i\otimes b_{p-i}=L_i\otimes M_{p-i}$ and Case \ref{caseLLYY} does
not occur\vspace*{-7pt}.
\end{case}

We claim that $L_i\otimes M_{p-i}-M_{p-i}\otimes L_i$ is the only
term in $r_p$ of the form $L_j\otimes M_{p-j}-M_{p-j}\otimes L_j$.
Indeed, let $i$ be maximal. Then it is not difficult to see that the
result holds for $i>0$.  If $i\leq0$, and suppose that there exists
$j_0\ne i$ such that $L_{j_0}\otimes M_{p-j_0}-M_{p-j_0}\otimes
L_{j_0}$ is a term in $r_p$ with nonzero coefficient. Take $j_0$ to
be minimal. Then $L_i\otimes [M_{p-i},L_{j_0}]\otimes M_{p-{j_0}}$
is a term in $\fc(r)_{2p}$ with nonzero coefficient, a
contradiction. This proves the claim. Now the condition
$\fc(r)_{2p}=0$ yields  that $\fc(L_i\otimes M_{p-i}-M_{p-i}\otimes
L_i)=0$. It follows that $i=0$ or $p-i=0$, i.e., $L_i\otimes
M_{p-i}-M_{p-i}\otimes L_i$ is equal to $L_p\otimes M_0-M_0\otimes
L_p$ or $L_0\otimes M_p-M_p\otimes L_0$. Now by our assumption that
Case 1 does not occur, we need only to consider terms of the form
$M_k\otimes M_{p-k}-M_{p-k}\otimes M_k$ for all $k\in\Z$.  It is not
difficult to check that the only possible term of the form
$M_k\otimes M_{p-k}-M_{p-k}\otimes M_k$ is $M_p\otimes
M_0-M_0\otimes M_p$. Thus $r_p\in V_3\cup V_4$\vspace*{-7pt}.

\begin{case}$a_i\otimes b_{p-i}=L_i\otimes Y_{p-i}$ or $M_i\otimes Y_{p-i}$\vspace*{-9pt}.\end{case}

In this case, $p\in\frac{1}{2}+\Z$ and $r_p\in{\rm Span}\{L_j\otimes
Y_{p-j}-Y_{p-j}\otimes L_j,M_j\otimes Y_{p-j}-Y_{p-j}\otimes
M_j\,|\,j\in\Z\}$. Assume that $a_i\otimes b_{p-i}=L_i\otimes
Y_{p-i}$ and  $i$ is the maximal integer such that $L_i\otimes
Y_{p-i}-Y_{p-i}\otimes L_i$ is a nonzero term in $r_p$. If $i>0$,
then we conclude that only $j\leq 0$ and $2(p-j)=i$ can $L_j\otimes
Y_{p-j}-Y_{p-j}\otimes L_j$ be a nonzero term in $r_p$. Meanwhile,
the condition $\fc(r)_{2p}=0$ yields the coefficient of $L_i\otimes
[Y_{p-i},Y_{p-j}]\otimes L_j$ in $\fc(r)_{2p}$ to be zero and so is
the coefficient of $L_j\otimes Y_{p-j}-Y_{p-j}\otimes L_j$. If
$i\leq 0$, then it is easy to see that the coefficient of
$L_j\otimes Y_{p-j}-Y_{p-j}\otimes L_j$ with $j\neq i$ must be zero.
Thus, we conclude that $L_i\otimes Y_{p-i}-Y_{p-i}\otimes L_i$ is
the only term of the form $L_k\otimes Y_{p-k}-Y_{p-k}\otimes L_k$ in
$r_p$. While for the case $a_i\otimes b_{p-i}=M_i\otimes Y_{p-i}$,
the same result can be obtained. Thus, $\fc(L_i\otimes
Y_{p-i}-Y_{p-i}\otimes L_i)=0$. It follows from that  $i=0$ or
$2(p-i)=i$, i.e., $L_i\otimes Y_{p-i}-Y_{p-i}\otimes L_i$ is equal
to $L_0\otimes Y_p-Y_p\otimes L_0$ or $L_i\otimes Y_{\frac
i2}-Y_{\frac i2}\otimes L_i$ with $i=\frac 2 3 p\in\Z$.  In the
former case, $M_0\otimes Y_p-Y_p\otimes M_0$ is the only possible
nonzero term of the form $M_k\otimes Y_{p-k}-Y_{p-k}\otimes
M_k(k\in\Z)$ in $r_p$, while in the latter case, the only
possibility is $M_i\otimes Y_{\frac i2}-Y_{\frac i2}\otimes M_i$.
Thus we conclude that $r_p\in V_6\cup V_7$ if $a_i\otimes
b_{p-i}=L_i\otimes Y_{p-i}$, otherwise $r_p\in V_8^{(i)}$
\vspace*{-7pt}.

\begin{case}$a_i\otimes b_{p-i}=M_i\otimes M_{p-i}$ and Cases \ref{caseLLYY} and
\ref{caseLM} do not occur\vspace*{-9pt}.\end{case}

Then one must have $r_p\in\sum_{j\in\Z}\F(M_j\otimes
M_{p-j}-M_{p-j}\otimes M_j)$, i.e., $r_p\in V_5$.\vspace*{4pt}

Now consider the case $p=0$.  By the similar argument as in Case
\ref{caseLLYY} one knows that the terms of the form
$Y_{i-\frac12}\otimes Y_{\frac12-i}-Y_{\frac12-i}\otimes
Y_{i-\frac12}(i\in\Z)$ cannot occur in $r_0$. So $r_0$ is in the
subspace spanned by
$$\{L_i\otimes L_{-i}-L_{-i}\otimes L_i,L_i\otimes
M_{-i}-M_{-i}\otimes L_i,M_i\otimes M_{-i}-M_{-i}\otimes
M_i\,|\,i\in \Z.\}$$  Let $i$ be the maximal index such that
$a_i\otimes b_{-i}-b_{-i}\otimes a_i$ is a term with nonzero
coefficient.   We may assume that $i>0$, since the case $i=0$ is
trivial. If $L_j\otimes L_{-j}-L_{-j}\otimes L_j$ is a term of $r_p$
with nonzero coefficient and $j\in\Z_{>0}$, then by the similar
argument as in Case 1 one has that $L_j\otimes L_{-j}-L_{-j}\otimes
L_j$  is the unique term of the form $L_k\otimes
L_{-k}-L_{-k}\otimes L_k$ with $k\in\Z_{>0}$ and $\fc(L_j\otimes
L_{-j}-L_{-j}\otimes L_j)=0$. But this implies $j=0$, contradicting
the choice of $j$. Similarly, for each $j\in\Z^*$ we deduce that
$L_j\otimes M_{-j}-M_{-j}\otimes L_j$ cannot be a term of $r_0$.
Thus we conclude that $r_0\in$ Span$\{L_0\otimes M_0-M_0\otimes L_0,
M_j\otimes M_{-j}-M_{-j}\otimes M_j\,|\,j\in\Z\}$. Now one can
easily see that either $r_0\in V_1$ or $r_0\in V_5$.\QED\vskip5pt

Although not all solutions to CYBE in $\LL$ can be solved (even in
the Witt algebra spanned by the set $\{L_i\,|\,i\in\Z\}$,
cf.~\cite{NT}),
 Proposition \ref{props} nevertheless provides us some rule to decide when
$r\in$ Im$(1-\tau)$ is not a solution to CYBE. Indeed, Proposition
\ref{props}  classifies all possible highest components $r_p$ of $r$
for which $r\in$ $\mbox{Im}(1-\tau)$ and $\fc(r)=0$. Similarly the
form of the lowest components $r_q$ can also be determined.

\vs{26pt}

\cl{\bf\S3. \ Proof of Theorem \ref{theo}}\setcounter{section}{3}
\setcounter{theo}{0}\setcounter{equation}{0}

\vs{8pt}

Regard $\VV=\LL\otimes\LL$ as a $\LL$-module under the adjoint
diagonal action. Denote by $\Der(\LL,\VV)$ the set of
\textit{derivations} $D:\LL\to\VV$, namely, $D$ is a linear map
satisfying
\begin{eqnarray}
\label{deriv} D([x,y])=x\cdot D(y)-y\cdot D(x),
\end{eqnarray}
and $\Inn(\LL,\VV)$ the set consisting of the derivations $v_{\rm
inn}, \, v\in\VV$, where $v_{\rm inn}$ is the \textit{inner
derivation} defined by $v_{\rm inn}:x\mapsto x\cdot v.$ Then it is
well known that $H^1(\LL,\VV)\cong\Der(\LL,\VV)/\Inn(\LL,\VV),
$
where $H^1(\LL,\VV)$ is the {\it first cohomology group} of the Lie
algebra $\LL$ with coefficients in the $\LL$-module $\VV$.

 A derivation
$D\in\Der(\LL,\VV)$ is {\it homogeneous of degree
$\alpha\in\frac12\Z$} if $D(\LL_p)\subset \VV_{\alpha+p}$ for all
$p\in\frac12\mathbb{Z}$.
 Denote $\Der(\LL,\VV)_\alpha =\{D\in\Der(\LL,\VV)\,|\,{\rm deg\,}D=
\alpha\}$ for $\alpha\in\frac{1}{2}\Z.$ Let $D$ be an element of
$\Der(\LL,\VV)$. For any $\alpha\in\frac12\Z$, define the linear map
$D_\alpha:\LL\rightarrow\VV$ as follows: For any $\mu\in\LL_q$ with
$q\in\frac12\mathbb{Z}$, write
$D(\mu)=\sum_{p\in\frac12\mathbb{Z}}\mu_p$ with $\mu_p\in\VV_p$,
then we set $D_\alpha(\mu)=\mu_{q+\alpha}$. Obviously, $D_\alpha\in
\Der(\LL,\VV)_\alpha$ and we have
\begin{eqnarray}\label{summable}
D=\mbox{$\sum\limits_{\alpha\in\frac12\mathbb{Z}}D_\alpha$},
\end{eqnarray}
which holds in the sense that for every $u\in\LL$, only finitely
many $D_\alpha(u)\neq 0,$ and
$D(u)=\sum_{\alpha\in\frac12\mathbb{Z}}D_\alpha(u)$ (we call such a
sum in (\ref{summable}) {\it summable}).

First we claim that if $\alpha\in\frac12\Z^*$ then
$D_\alpha\in\Inn(\LL,\VV)$. To see this, denote
$\gamma=\alpha^{-1}D_{\alpha}(L_0)\in\VV_{\alpha}$. Then for any
$x_n\in\LL_{n}$, applying $D_{\alpha}$ to $[L_0,x_n]=nx_n$ and using
$D_{\alpha}(x_n)\in \VV_{n+\alpha}$, we obtain
$(\alpha+n)D_{\alpha}(x_n)-x_n\cdot D_{\alpha}(L_0)=L_0\cdot
D_{\alpha}(x_n)-x_n\cdot D_{\alpha}(L_0)=nD_{\alpha}(x_n),
$
 i.e.,
$D_{\alpha}(x_n)=\gamma_{\rm inn}(x_n)$. Thus
$D_{\alpha}=\gamma_{\rm inn}$ is inner.

 In the following we always
use the symbol  ``$\equiv$''  to denote modulo $\F(M_0\otimes M_0)$.
Then we can claim that $D_0(L_0)\equiv0$. Indeed, for any
$p\in\frac12\Z$ and $x_p\in\LL_p$, applying $D_0$ to
$[L_0,x_p]=px_p$, one has $x_p\cdot D_0(L_0)=0$. Thus by Lemma
\ref{Legr}, $D_0(L_0)\equiv0$.

Now we claim that for any $D\in{\rm Der}(\LL,\VV)$, (\ref{summable})
is a finite sum. To see this, one can suppose $D_n=(v_n)_{\rm inn}$
for some $v_n\in\VV_n$ and $n\in\frac{1}{2}\Z^*$. If
$\Z'=\{n\in\frac{1}{2}\Z^*\,|\,v_n\ne0\}$ is an infinite set, then
$D(L_0)=D_0(L_0)+\sum_{n\in\Z'}L_0\cdot
v_n=D_0(L_0)+\sum_{n\in\Z'}nv_n$ is an infinite sum, which is not an
element in $\VV$, contradicting
 the fact that $D$ is a derivation from $\LL$ to $\VV$. This
 together with the proposition below
proves  Theorem \ref{theo}(i). (To complete proof of Theorem
\ref{theo}(i), one still needs to show ${\rm
Inn}(\LL,\VV)\cap\DD=\{0\}$. For this, suppose $D=u_{\rm
inn}\in\DD$, where $u$ is a linear combination of $a_i\otimes b_j$
for some $a_i,b_j\in\{L_m,Y_{m-\frac12},M_m\,|\,m\in\Z\}$. Applying
$D$ to generators of $\LL$ and using (\ref{def-D}), one immediately
obtains $D=0$.)
\begin{prop}\label{p1}
\rm Replacing $D_0$ by $D_0-u_{\rm inn}$ for some $u\in \VV_0$, one
can suppose $D_0\in\DD$ $($where $\DD$ is defined in
$(\ref{def-D}))$.
\end{prop}
\noindent{\it Proof}\ \ The proof seems to be technical. We shall
prove that after a number of steps in each of which $D_0$ is
replaced by $D_0-u_{\rm inn}$ for some $u\in \VV_0$, we obtain
$D_0\in \DD$. This will be done by some lengthy calculations.

For any $n\in\Z$, one can write $D_0(Y_{n-\frac12})$, $D_0(L_n)$ and
$D_0(M_n)$ as follows
\begin{eqnarray*}
&&\!\!\!\!\!\!\!\!\!\!\!\!
D_0(Y_{n-\frac12})\!=\!\mbox{$\sum\limits_{i\in\Z}$}(\a_{n,i}L_i\!\otimes\!
Y_{n-\frac12-i}\!+\!\a^\dag_{n,i}Y_{i-\frac12}\!\otimes\!
L_{n-i}\!+\!\b_{n,i}M_i\!\otimes\!
Y_{n-\frac12-i}\!+\!\b^\dag_{n,i}Y_{i-\frac12}\!\otimes\!M_{n-i}),\\
&&\!\!\!\!\!\!\!\!\!\!\!\!
D_0(L_n)\!=\!\mbox{$\sum\limits_{i\in\Z}$}(a_{n,i}L_i\!\otimes\!
L_{n-i}\!+\!b_{n,i}L_i\!\otimes\!M_{n-i}\!+\!b^\dag_{n,i}M_i\!\otimes\!
L_{n-i}\!+\!c_{n,i}M_i\!\otimes\!M_{n-i}\!+\!d_{n,i}Y_{i-\frac{1}{2}}\!\otimes\!Y_{n-i+\frac{1}{2}}),\\
&&\!\!\!\!\!\!\!\!\!\!\!\!
D_0(M_n)\!=\!\mbox{$\sum\limits_{i\in\Z}$}(e_{n,i}L_i\!\otimes\!
L_{n-i}\!+\!f_{n,i}L_i\!\otimes\!M_{n-i}\!+\!f^\dag_{n,i}M_i\!\otimes\!
L_{n-i}\!+\!g_{n,i}M_i\!\otimes\!
M_{n-i}\!+\!h_{n,i}Y_{i-\frac{1}{2}}\!\otimes\!Y_{n-i+\frac{1}{2}}),
\end{eqnarray*}
where all  coefficients of the tensor products are in $\F$, and the
sums are all finite. For any $n\in\Z$, the following identities
hold,
\begin{eqnarray*}
&&L_1\cdot(M_n\otimes M_{-n})=nM_{n+1}\otimes
M_{-n}-nM_n\otimes M_{1-n},\\
&&L_1\cdot(L_n\otimes M_{-n})=(n-1)L_{n+1}\otimes
M_{-n}-nL_n\otimes M_{1-n},\\
&&L_1\cdot(M_n\otimes L_{-n})=nM_{n+1}\otimes
L_{-n}-(1+n)M_n\otimes L_{1-n},\\
&&L_1\cdot(L_n\otimes L_{-n})=(n-1)L_{n+1}\otimes
L_{-n}-(1+n)L_n\otimes L_{1-n},\\
&&L_1\cdot(Y_{n-\frac12}\otimes
Y_{\frac12-n})=(n-1)Y_{n+\frac12}\otimes
 Y_{\frac12-n}-nY_{n-\frac12}\otimes Y_{3/2-n}.
\end{eqnarray*}
Let $\triangle$ denote the set consisting of 5 symbols
$a,b,b^\dag,c,d$. For each $x\in \triangle$ we define
$M_{x}=\max\{\,|p\,|\,\big|\,x_{1,p}\ne0\}.$ For $n=1$, using the
induction on $\sum_{x\in \triangle} M_x$ in the above identities,
and replacing $D_0$ by $D_0-u_{\rm inn}$, where $u$ is a proper
linear combination of $L_{p}\otimes L_{-p}$, $L_{p}\otimes M_{-p}$,
$M_{p}\otimes L_{-p}$, $M_{p}\otimes M_{-p}$ and
$Y_{p-\frac12}\otimes Y_{\frac12-p}$ with $p\in\Z$, one can safely
suppose
\begin{eqnarray}\label{assumption}
a_{1,i}=b_{1,j}=b^\dag_{1,k}=c_{1,m}=d_{1,n}=0\ {\rm for}\
i\neq-1,2,\,j\neq 0,2,\,k\neq\pm1,\ m\neq 0,1,\,n\neq 0,2.
\end{eqnarray}
Applying $D_0$ to $[\,L_{1},L_{-1}]=-2L_0$ and using the fact that
$D_0(L_0)=dM_0\otimes M_0$ for some $d\in \F$, we obtain
\vspace{-5pt}
\begin{eqnarray*}\label{L1L-1}
\begin{aligned}
&\mbox{$\sum\limits_{p\in
\Z}$}\Big(\big((p-2)a_{-1,p-1}-(p+2)a_{-1,p}+(p-2)a_{1,p}-(p+2)a_{1,1+p}\big)L_p\otimes
L_{-p}\\
&+\big((p-2)b_{-1,p-1}-(1+p)b_{-1,p}+(p-1)b_{1,p}-(p+2)b_{1,1+p}\big)L_p\otimes
M_{-p}\\
&+\big((p-1)b^\dag_{-1,p-1}-(p+2)b^\dag_{-1,p}
+(p-2)b^\dag_{1,p}-(p+1)b^\dag_{1,1+p}\big)M_p\otimes L_{-p}\\
&+\big((p-1)c_{-1,p-1}-(p+1)c_{-1,p}
+(p-1)c_{1,p}-(p+1)c_{1,p+1}+2\delta_{0,p}d\big)M_p\otimes M_{-p}\\
&+\big((p-2)d_{-1,p-1}-(1+p)d_{-1,p}
+(p-2)d_{1,p}-(p+1)d_{1,1+p}\big)Y_{p-\frac12}\otimes
Y_{\frac12-p}\Big)=0.
\end{aligned}
\end{eqnarray*}
In particular, for any $p\in\Z$ one has
\begin{eqnarray*}
(p-2)a_{-1,p-1}-(p+2)a_{-1,p}+(p-2)a_{1,p}-(p+2)a_{1,1+p}=0,
\end{eqnarray*}
which together with the fact that $\{p\in \Z\,|\,a_{-1,p}\neq0\}$ is
finite, forces
\begin{eqnarray}\label{a-1}
a_{-1,p}=a_{-1,0}+3a_{-1,1}+3a_{1,2} =3a_{-1,-2}+a_{-1,-1}+3a_{1,-1}
=a_{-1,-1}+a_{-1,0}=0
\end{eqnarray}
for  $p\in\Z\SM\{-2,0,\pm1\}$. Similarly, comparing the coefficients
of $L_p\otimes M_{-p}$, $M_p\otimes L_{-p}$ , $M_p\otimes M_{-p}$
and $Y_{p-\frac12}\otimes Y_{\frac12-p}$, one has
\begin{eqnarray}
\!\!\!\!\!\!\!\!\!\!\!\!\!\!\!&&\!\!\!\!\!b_{1,q}=b^\dag_{1,q}=b_{-1,p_1}=b^\dag_{-1,p_2}=c_{-1,p_3}=d_{-1,p_4}=0,\nonumber\\
\!\!\!\!\!\!\!\!\!\!\!\!\!\!\!&&\!\!\!\!\!\!\!\!\begin{array}{lll}\label{bcd-1j}
b_{-1,0}\!+\!2b_{-1,-1}=b_{-1,1}\!-\!b_{-1,-1}=b^\dag_{-1,-1}\!+\!2b^\dag_{-1,-2}
=c_{-1,-1}\!+\!c_{-1,0}\!+\!c_{1,1}\!+\!c_{1,0}-2d\\[4pt]
\ \ \ \ \ \ \ \ \ \ \,\ \ \ \ \ \ \
\,\,=b^\dag_{-1,0}\!-\!b^\dag_{-1,-2}=2d_{-1,-1}\!+\!d_{-1,0}\!+\!2d_{1,0}
=d_{-1,0}\!+\!2d_{-1,1}\!+\!2d_{1,2}=0,
\end{array}
\end{eqnarray}
for any $q\in\Z$, $p_1\in\Z\SM\{0,\pm1\}$, $p_2\in\Z\SM\{0,-1,-2\}$,
$\forall\ p_3\in\Z\SM\{-1,0\}$ and $p_4\in\Z\SM\{0,\pm 1\}$.

Note that $L_1\cdot(L_1\otimes M_{-1}-L_0\otimes M_0)=0$. Replacing
$D_0$ by $D_0+b_{-1,1}(L_1\otimes M_{-1}-L_0\otimes M_0)$, one can
assume $b_{-1,-1}=b_{-1,0}=b_{-1,1}=0$ by (\ref{bcd-1j}). Similarly,
replacing $D_0$ by $D_0+b_{-1,-2}^\dag(M_{-1}\otimes L_1-M_0\otimes
L_0)$, one can suppose
$b_{-1,-2}^\dag=b_{-1,-1}^\dag=b_{-1,0}^\dag=0$. Hence
$b_{-1,p}=b_{1,p}=b^\dag_{-1,p}=b^\dag_{1,p}=0$ for all $p\in\Z.$
Applying $D_0$ to $[\,L_{2},L_{-1}]=-3L_{1}$, one has

\begin{eqnarray*}\label{L2L-1}
\begin{aligned}
&\mbox{$\sum\limits_{p\in
\Z}$}\Big((p-4)a_{-1,p-2}-(3+p)a_{-1,p}-(p+2)a_{2,p+1}-(3-p)a_{2,p}+3a_{1,p}\big)L_{p}\otimes
L_{1-p}\\
&-\big((p+2)b_{2,p+1}+(2-p)b_{2,p}\big)L_{p}\otimes M_{1-p}-\big((p+1)b^\dag_{2,p+1}+(3-p)b^\dag_{2,p}\big)M_{p}\otimes L_{1-p}\\
&+\big((p-2)c_{-1,p-2}-(1+p)c_{-1,p}-(p+1)c_{2,p+1}-(2-p)c_{2,p}+3c_{1,p}\big)M_{p}\otimes M_{1-p}\\
&+\big((p-7/2)d_{-1,p-2}\!-\!(3/2+p)d_{-1,p}\!-\!(p+1)d_{2,p+1}\!-\!(3-p)d_{2,p}\!+\!3d_{1,p}\big)Y_{p-\frac12}\otimes
Y_{3/2-p}\Big)=0.
\end{aligned}
\end{eqnarray*}
 By  computing the
coefficient of $L_p\otimes L_{1-p}$ and using (\ref{assumption}),
(\ref{a-1}), and that $\{p\,|\,a_{2,p}\neq0\}$ is finite,  one has
\begin{eqnarray}\label{a2}
0\!\!\!\!\!\!\!&&=a_{2,p}=a_{-1,-2}=a_{-1,1}\ \ \  \mathrm{for}\ p\in\Z\SM\{0,\pm1,2,3\}\nonumber\\
&&=4a_{2,-1}-(3a_{-1,0}-a_{2,0})=2a_{2,1}+3(a_{-1,0}+a_{2,0})\\
&&=a_{2,2}-(2a_{-1,0}+a_{2,0})=4a_{2,3}+(5a_{-1,0}+a_{2,0}).\nonumber
\end{eqnarray}
Similarly, one can obtain that
\begin{eqnarray}
\label{d-1}&d_{-1,-1}=d_{-1,1}=0,\\
\label{bcd2}&b_{2,p_1}=b^\dag_{2,p_2}=c_{2,p_3}=d_{2,p_4}=0,
\end{eqnarray}
\vskip-20pt\begin{eqnarray}\label{bcd2j}
\begin{aligned}
&b_{2,0}+3b_{2,-1}=b_{2,1}-3b_{2,-1}=b_{2,2}+b_{2,-1}=
b^\dag_{2,0}+b^\dag_{2,3}=b^\dag_{2,1}-3b^\dag_{2,3}=b^\dag_{2,2}+3b^\dag_{2,3}\\
&\ \ \ \ \ \ \ \ \ \ \ \ \ \ \ =c_{2,1}-(3c_{1,0}-c_{-1,0}-2c_{2,0})\\
&\ \ \ \ \ \ \ \ \ \ \ \ \ \ \ =2c_{2,2}-(3c_{1,1}-c_{-1,-1}-3c_{1,0}+c_{-1,0}+2c_{2,0})\\
&\ \ \ \ \ \ \ \ \ \ \ \ \ \ \  =d_{2,1}-3(2d_{1,0}-d_{2,0})=
d_{2,2}+3(2d_{1,0}-d_{2,0})=d_{2,3}-(4d_{1,0}-d_{2,0})=0,
\end{aligned}
\end{eqnarray}
for any $ p_1\in\Z\SM\{\pm1,0,2\},
p_2\in\Z\SM\{0,1,2,3\},p_3\in\Z\SM\{0,1,2\}$ and
$p_4\in\Z\SM\{0,1,2,3\}$. From the equation
$[\,L_{1},L_{-2}]=-3L_{-1}$, we obtain
\begin{eqnarray*}\label{L1L-2}
\begin{aligned}
 &\mbox{$\sum\limits_{p\in
\Z}$}\Big((p-2)a_{-2,p-1}-(3+p)a_{-2,p}-(p+4)a_{1,p+2}+(p-3)a_{1,p}+3a_{-1,p}\big)L_{p}\otimes
L_{-1-p}\\
&+\big((p-2)b_{-2,p-1}-(2+p)b_{-2,p}\big)L_{p}\otimes
M_{-1-p}+\big((p-1)b^\dag_{-2,p-1}-(3+p)b^\dag_{-2,p}\big)M_{p}\otimes
L_{-1-p}\\
&+\big((p-1)c_{-2,p-1}\!-\!(2+p)c_{-2,p}-(p+2)c_{1,p+2}+(p-1)c_{1,p}+3c_{-1,p}\big)M_{p}\otimes M_{-1-p}\\
&+\big((p\!-\!2)d_{-2,p-1}\!-\!(2\!+\!p)d_{-2,p}\!-\!(p+5/2)d_{1,p+2}\!+\!(p-5/2)d_{1,p}\!+\!3d_{-1,p}\big)Y_{p-\frac12}\!\otimes\!
Y_{-p-\frac12}\Big)=0.
\end{aligned}
\end{eqnarray*}
It follows from the above formula  and
(\ref{assumption})--(\ref{d-1}) that
\begin{eqnarray}
\label{a+-1=0}&a_{\pm1,p}=d_{\pm1,p}=0,\\
\label{abcd-2}&a_{-2,p_1}=b_{-2,p_2}=b^\dag_{-2,p_3}=c_{-2,p_4}=d_{-2,p_5}=0,
\end{eqnarray}
\vskip-20pt
\begin{eqnarray}\label{abcd-2j}
\begin{aligned}
&a_{-2,-3}-a_{-2,1}=a_{-2,-2}+4a_{-2,1}=a_{-2,-1}-6a_{-2,1}=a_{-2,0}+4a_{-2,1}\\
&\ \ \ \ \ \ \ \ \ \ \ \ \ \ \ \ \  =b_{-2,-1}+3b_{-2,-2}=b_{-2,0}-3b_{-2,-2}=b_{-2,1}+b_{-2,-2}\\
&\ \ \ \ \ \ \ \ \ \ \ \ \ \ \ \ \  =b^\dag_{-2,-2}+3b^\dag_{-2,-3}=b^\dag_{-2,-1}-3b^\dag_{-2,-3}=b^\dag_{-2,0}+b^\dag_{-2,-3}\\
&\ \ \ \ \ \ \ \ \ \ \ \ \ \ \ \ \  =c_{-2,-1}-(3c_{-1,-1}-2c_{-2,-2}-c_{1,1})\\
&\ \ \ \ \ \ \ \ \ \ \ \ \ \ \ \ \  =2c_{-2,0}-(3c_{-1,0}-c_{1,0}-3c_{-1,-1}+2c_{-2,-2}+c_{1,1})\\
&\ \ \ \ \ \ \ \ \ \ \ \ \ \ \ \ \
=d_{-2,-1}+3d_{-2,-2}=d_{-2,0}-3d_{-2,-2}=d_{-2,1}+d_{-2,-2}=0,
\end{aligned}
\end{eqnarray}
for all $p\in \Z,\,
p_1\in\Z\SM\{-3,-2,0,\pm1\},\,p_2\in\Z\SM\{-2,0,\pm1\},\,p_3\in\Z\SM\{-3,-2,-1,0\},$
$p_4\in\Z\SM\{-2,-1,0\},\,p_5\in\Z\SM\{-2,0,\pm1\}$. Applying $D_0$
to $[L_2,L_{-2}]=-4L_0$, one has
\begin{eqnarray*}
\begin{aligned}
&\mbox{$\sum\limits_{p\in\Z}$}\Big(\big((p-4)a_{-2,p-2}-(p+4)a_{-2,p}-(p+4)a_{2,p+2}+(p-4)a_{2,p}\big)L_{p}\otimes L_{-p}\\
&+\big((p-4)b_{-2,p-2}-(p+2)b_{-2,p}-(p+4)b_{2,p+2}+(p-2)b_{2,p}\big)L_{p}\otimes M_{-p}\\
&+\big((p-2)b^\dag_{-2,p-2}\!-\!(p+4)b^\dag_{-2,p}\!-\!(p+2)b^\dag_{2,p+2}+(p-4)b^\dag_{2,p}\big)M_{p}\otimes L_{-p}\\
&+\big((p-2)c_{-2,p-2}\!-\!(p+2)c_{-2,p}-(p+2)c_{2,p+2}+(p-2)c_{2,p}+4d\delta_{0p}\big)M_{p}\otimes M_{-p}\\
&+\big((p-7/2)d_{-2,p-2}\!-\!(p+5/2)d_{-2,p}-(p+5/2)d_{2,p+2}+(p-7/2)d_{2,p}\big)d_{p-\frac12}\otimes
Y_{\frac12-p}\Big)=0,
\end{aligned}
\end{eqnarray*}
which combined with (\ref{a2}) and (\ref{bcd2})--(\ref{abcd-2j})
yields the follows:
\begin{eqnarray}
&b_{\pm2,p}=b^\dag_{\pm2,p}=0,\ \forall\,\,p\in\Z,\nonumber\\
\label{xin}&a_{-2,1}+a_{2,-1}=c_{-2,-1}+c_{2,1}=d_{-2,-2}+d_{2,0}=0,\\
&c_{2,0}+c_{2,2}+c_{-2,0}+c_{-2,-2}=2d\nonumber.
\end{eqnarray}

Set $u=L_{-1}\otimes L_1-2L_0\otimes L_0+L_1\otimes L_{-1}$. Observe
that $L_{\pm1}\cdot u=0,$ so, by equations (\ref{a2}),
(\ref{a+-1=0}), (\ref{abcd-2j}) and (\ref{xin}), one can assume
$$a_{2,-1}=a_{2,0}=a_{2,1}=a_{2,2}=a_{2,3}=a_{-2,-3}=a_{-2,-2}=a_{-2,-1}=a_{-2,0}=a_{-2,1}=0,$$
when $D_0$ is replaced by $D_0+a_{2,-1}(L_{-1}\otimes
L_1-2L_0\otimes L_0+L_1\otimes L_{-1}).$ Similarly, set
$u=Y_{\frac12}\otimes Y_{-\frac12}-Y_{-\frac12}\otimes Y_{\frac12}$,
one can assume
$d_{2,0}=d_{2,1}=d_{2,2}=d_{2,3}=d_{-2,-2}=d_{-2,-1}=d_{-2,0}=d_{-2,1}=0.$
Hence so far we have obtained that
$a_{\pm1,p}=a_{\pm2,p}=b_{\pm1,p}=b_{\pm2,p}=b^\dag_{\pm1,p}=b^\dag_{\pm2,p}=d_{\pm1,p}=d_{\pm2,p}
=0\ \mathrm{for\ any}\ p \in\Z$ and

\begin{eqnarray}
&&c_{-1,-1}+c_{-1,0}+c_{1,1}+c_{1,0}-2d=c_{1,p} = c_{-1,-p}=0  \ \ \forall\ p\in\Z\SM\{0,1\}\nonumber\\
\label{cz}&&c_{2,p}=c_{-2,-p}=0\ \  \  \forall\ p\in\Z\SM\{0,1,2\}\nonumber\\
&&\ \ \ \ =c_{2,1}-(3c_{1,0}-c_{-1,0}-2c_{2,0})=2c_{2,2}-(3c_{1,1}-c_{-1,-1}-3c_{1,0}+c_{-1,0}+2c_{2,0})\nonumber\\
&&\ \ \ \ =c_{-2,-1}-(3c_{-1,-1}-2c_{-2,-2}-c_{1,1})=2c_{-2,0}-(3c_{-1,0}-c_{1,0}-3c_{-1,-1}+2c_{-2,-2}+c_{1,1})\nonumber\\
&&\ \ \ \ =c_{-2,-1}+c_{2,1}=c_{2,0}+c_{2,2}+c_{-2,0}+c_{-2,-2}-2d.
\end{eqnarray}
It follows from by  repeatedly applying ad$L_1$ to $L_2$ and using
the fact that for all $n\in\Z$ $D_0(L_n)$ is a finite sum that
$c_{2,1}=0$. Hence by (\ref{cz}),  $c_{-2,-1}=0$. Now replacing
$c_{\pm1,0}$, $c_{\pm1,\pm1}$, $c_{\pm2,0}$ and $c_{\pm2,\pm2}$  by
$\gamma\pm\a$, $\gamma^\dag\pm\a^\dag$, $\gamma\pm 2\a$ and
$\gamma^\dag\pm 2\a^\dag$ in (\ref{cz}) respectively,
$D_0(L_{\pm1})$ and $D_0(L_{\pm2})$ have the following more concise
expressions:
\begin{eqnarray}
&D_0(L_\pm)=(\gamma\pm\a)M_0\otimes
M_{\pm1}+(\gamma^\dag\pm\a^\dag)M_{\pm1}\otimes M_0,\label{L1}\\
&D_0(L_{\pm2})=(\gamma\pm2\a)M_0\otimes
M_{\pm2}+(\gamma\pm2\a^\dag)M_{\pm2}\otimes M_0.\label{L2}
\end{eqnarray}
Thus for any $n\in \Z$ one can deduce
$D_0(L_n)=(n\a+\gamma)M_0\otimes
M_n+(n\a^\dag+\gamma^\dag))M_n\otimes M_0$, since $\WW$ can be
generated by $L_{\pm1}$ and $L_{\pm2}$.

To prove the proposition we still need to show $D_0(Y_{\frac12})=\b
M_0\otimes Y_{\frac12}+\b^\dag Y_{\frac12}\otimes M_0$. Applying
$D_0$ to $[L_m,Y_{n-\frac12}]=(n-(m+1)/2)Y_{n+m-\frac12}$ and
noticing that $Y_{n-\frac12}\cdot D_0L_m=0$, we obtain
\begin{eqnarray}
\begin{aligned}
&(i-2m)\a_{n,i-m}+(n-i-(m+1)/2)\a_{n,i}-(n-(m+1)/2)\a_{n+m,i}=0   ,\\
&(i-(3m+1)/2)\a^\dag_{n,i-m}+(n-m-i)\a^\dag_{n,i}-(n-(m+1)/2)\a^\dag_{n+m,i}=0 ,\\
&(i-m)\b_{n,i-m}+(n-i-(m+1)/2)\b_{n,i}-(n-(m+1)/2)\b_{n+m,i}=0 ,\\
&(i-(3m+1)/2)\b^\dag_{n,i-m}+(n-i)\b^\dag_{n,i}-(n-(m+1)/2)\b^\dag_{n+m,i}=0.
\end{aligned}
\end{eqnarray}
In the above equations, putting $n=m=1$ and using the fact that the
rank of $\{x_{1,p}\ |\ x=\a,\a^\dag,\b\ \mathrm{ or}\ \b^\dag\}$ is
finite, one has
\begin{eqnarray}
\begin{aligned}
&\a_{1,1}+\a_{1,0}=\a^\dag_{1,1}+\a^\dag_{1,0}=0,\\
&\a_{1,p_1}=\a^\dag_{1,p_1}=\b_{1,p_2}=\b^\dag_{1,p_3}=0,\ \ \forall\,\,p_1\in\Z\SM\{0,1\},\,p_2\in\Z\SM\{0\},\,p_3\in\Z\SM\{1\}.\\
\end{aligned}
\end{eqnarray}
 Similarly, letting $m=-1$ and $n=0$, then one has
\begin{eqnarray}
\begin{aligned}
&\a_{0,-1}+\a_{0,0}=\a^\dag_{0,1}+\a^\dag_{0,0}=0,\\
&\a_{0,p_1}=\a^\dag_{0,p_2}=\b_{0,p_3}=\b^\dag_{0,p_3}=0,
\ \ \forall\,\,p_1\in\Z\SM\{0,-1\},\,p_2\in\Z\SM\{0,1\},\,p_3\in\Z\SM\{0\}.\\
\end{aligned}
\end{eqnarray}
Taking $n=1$ and $m=-1$, one has $\a_{1,0}=-\a_{0,0}$,
$\a^\dag_{1,0}=\a^\dag_{0,0}$, $\b_{1,0}=\b_{0,0}$,
$\b^\dag_{1,1}=\b^\dag_{0,0}$. Thus $D_0(Y_{\frac12})$ and
$D_0(Y_{-\frac12})$ can be written as
\begin{eqnarray*}
D_0(Y_{\frac12})\!\!\!&=&\!\!\!-\a_{0,0}L_0\otimes
Y_{\frac12}+\a_{0,0}L_1\otimes
Y_{-\frac12}+\a^\dag_{0,0}Y_{-\frac12}\otimes
L_{1}-\a^\dag_{0,0}Y_{\frac12}\otimes L_{0}\\
\!\!\!&&\!\!\!+\b_{0,0}M_0\otimes
Y_{\frac12}+\b^\dag_{0,0}Y_{\frac12}\otimes M_0,\\
D_0(Y_{-\frac12})\!\!\!&=&\!\!\!-\a_{0,0}L_{-1}\otimes
Y_{\frac12}+\a_{0,0}L_0\otimes
Y_{-\frac12}-\a^\dag_{0,0}Y_{\frac12}\otimes
L_{-1}+\a^\dag_{0,0}Y_{-\frac12}\otimes L_{0}\\
\!\!\!&&\!\!\!+ \b_{0,0}M_0\otimes
Y_{-\frac12}+\b^\dag_{0,0}Y_{-\frac12}\otimes M_0.
\end{eqnarray*}
While for $M_n$ we have
\begin{eqnarray}
\begin{aligned}
&(i\!-\!2m)e_{n,i-m}\!+\!(n\!-\!m\!-\!i)e_{n,i}\!-\!ne_{n+m,i}=(i-2m)f_{n,i-m}\!+\!(n-i)f_{n,i}-nf_{n+m,i}=0   ,\\
&(i-(3m+1)/2)h_{n,i-m}+(n-i-(m-1)/2)h_{n,i}-nh_{n+m,i}=0,\\
&(i-m)f^\dag_{n,i-m}\!+\!(n-m-i)f^\dag_{n,i}\!-nf^\dag_{n+m,i}=(i-m)g_{n,i-m}\!+\!(n-i)g_{n,i}-ng_{n+m,i}=0.
\end{aligned}
\end{eqnarray}
For fixed $n=0$, putting $m=1$ and $m=-1$, respectively, one can
deduce
\begin{eqnarray*}
&&e_{0,-1}-e_{0,1}\!=\!e_{0,0}+2e_{0,1}\!=\!h_{0,1}+h_{0,0}\!=\!0,\\
&&e_{0,p_1}\!=\!f_{0,p_2}=f^\dag_{0,p_2}\!=\!g_{0,p_3}\!=\!h_{0,p_4}\!=\!0,\,\,
\forall\,\,p_1\in\Z\SM\{0,\pm1\},\,p_2\in\Z,\,p_3\in\Z^*,\,p_4\in\Z\SM\{0,1\}.
\end{eqnarray*}
Thus $D_0(M_0)$ can be written as
\begin{eqnarray*}
D_0(M_0)\!\!\!&=&\!\!\!e_{0,1}L_{-1}\otimes L_{1}
-2e_{0,1}L_{0}\otimes L_{0}+e_{0,1}L_1\otimes L_{-1}
+g_{0,0}M_0\otimes M_0\\
\!\!\!&&\!\!\!-h_{0,1}Y_{-\frac12}\otimes
Y_{\frac12}+h_{0,1}Y_{\frac12}\otimes Y_{-\frac12}.
\end{eqnarray*}
Applying $D_0$ to the equation $[Y_{-\frac12},Y_{\frac12}]=M_0$, we
get
\begin{eqnarray*}
e_{0,1}=0,g_{0,0}=2(\b_{0,0}+\b^\dag_{0,0}),h_{0,1}=\frac32(\a_{0,0}-\a^\dag_{0,0}).
\end{eqnarray*}
Thus, one can rewrite $D_0(M_0)$ as
\begin{eqnarray*}
D_0(M_0)=2(\b_{0,0}+\b^\dag_{0,0})M_0\otimes
M_0-h_{0,1}Y_{-\frac12}\otimes Y_{\frac12}+h_{0,1}Y_{\frac12}\otimes
Y_{-\frac12}.
\end{eqnarray*}
Notice that $[M_0,Y_{\frac12}]=0$, applying $D_0$ to which one would
have $h_{0,1}=0$, that is to say $\a_{0,0}=\a^\dag_{0,0}$. So we can
further simplify $D_0(Y_{\frac12})$, $D_0(Y_{-\frac12})$ and
$D_0(M_0)$ as follows:
\begin{eqnarray*}
D_0(M_0)\!\!\!&=&\!\!\!2(\b_{0,0}+\b^\dag_{0,0})M_0\otimes M_0,\\
D_0(Y_{\frac12})\!\!\!&=&\!\!\!-\a_{0,0}L_0\otimes
Y_{\frac12}+\a_{0,0}L_1\otimes
Y_{-\frac12}+\a_{0,0}Y_{-\frac12}\otimes L_{1}\\
\!\!\!&&\!\!\!-\a_{0,0}Y_{\frac12}\otimes L_{0}+ \b_{0,0}M_0\otimes
Y_{\frac12}+\b^\dag_{0,0}Y_{\frac12}\otimes M_0,\\
D_0(Y_{-\frac12})\!\!\!&=&\!\!\!-\a_{0,0}L_{-1}\otimes
Y_{\frac12}+\a_{0,0}L_0\otimes
Y_{-\frac12}-\a_{0,0}Y_{\frac12}\otimes L_{-1}\\
\!\!\!&&\!\!\!+\a_{0,0}Y_{-\frac12}\otimes L_{0}+ \b_{0,0}M_0\otimes
Y_{-\frac12}+\b^\dag_{0,0}Y_{-\frac12}\otimes M_0.
\end{eqnarray*}
Using the equation $[L_2,Y_{-\frac12}]=-\frac32Y_{3/2}$, we can
deduce
\begin{eqnarray*}
D_0(Y_{3/2})\!\!\!&=&\!\!\!-2\a_{0,0}L_1\otimes
Y_{\frac12}-2\a_{0,0}Y_{\frac12}\otimes
L_1-\frac{\a_{0,0}}3L_{-1}\otimes Y_{5/2}
-\frac{\a_{0,0}}3 Y_{5/2}\otimes L_{-1}\\
\!\!\!&&\!\!\!+\frac{4\a_{0,0}}3L_2\otimes
Y_{-\frac12}+\frac{4\a_{0,0}}3Y_{-\frac12}\otimes L_2+
\a_{0,0}L_0\otimes Y_{3/2}+
\a_{0,0}Y_{3/2}\otimes L_0\\
\!\!\!&&\!\!\!+ \b_{0,0}M_0\otimes Y_{3/2}+
\b^\dag_{0,0}Y_{3/2}\otimes M_0.
\end{eqnarray*}
Applying $D_0$ to $[L_{-2},Y_{3/2}]=\frac52Y_{-\frac12}$ and
noticing $Y_{3/2}\cdot D_0(L_{-2})=0$, one has $\a_{0,0}=0$, which
yields
\begin{eqnarray*}
D_0(Y_{\frac12})\!\!\!&=&\!\!\!\b_{0,0}M_0\otimes
Y_{\frac12}+\b^\dag_{0,0}Y_{\frac12}\otimes M_0.
\end{eqnarray*}
Now the  statement in Proposition \ref{p1} can be obtained
immediately, since $L_{\pm1}$, $L_{\pm2}$ and $Y_{\frac12}$ is  a
system of generators of $\LL$.\QED\vskip5pt

To prove the second part of the main theorem, we need the following
lemma.
\begin{lemm}\rm \label{lemma3ll}
Suppose $v\in\VV$ such that $x\cdot v\in {\rm Im}(1-\tau)$ for all
$x\in\LL.$ Then $v-d_0M_0\otimes M_0\in {\rm Im}(1-\tau)$ for some
$d_0\in\F$.
\end{lemm}
\ni{\it Proof}\ \ First note that $\LL\cdot {\rm Im}(1-\tau)\subset
{\rm Im}(1-\tau).$  We prove that after several steps, by replacing
$v$ with $v-u$ for some $u\in {\rm Im}(1-\tau)$, we shall have
$v-d_0M_0\otimes M_0=0$ for some $d_0\in \F$ and thus $v-d_0
M_0\otimes M_0\in{\rm Im}(1-\tau)$. Write
$v=\sum_{n\in\frac{1}{2}\Z}v_n.$ Obviously,
\begin{eqnarray}\label{eqrx}
v\in {\rm Im}(1-\tau)\ \,\Longleftrightarrow \ \,v_n\in {\rm
Im}(1-\tau),\ \ \forall\,\,n\in\frac{1}{2}\Z.
\end{eqnarray}
Then $\sum_{n\in\frac{1}{2}\Z}nv_n=L_0\cdot v\in {\rm Im}(1-\tau)$.
By (\ref{eqrx}), $nv_n\in {\rm Im}(1-\tau),$ in particular,
$v_{n}\in {\rm Im}(1-\tau)$ if $n\ne0$. Thus by replacing $v$ by
$v-\sum_{n\in\frac{1}{2}\Z^*}v_n$, one can suppose $v=v_0\in\VV_0$.
Write
\begin{eqnarray*}
v\!\!\!&=&\!\!\!\mbox{$\sum\limits_{p\in\Z}$}(a_{p}L_p\otimes
L_{-p}+b_{p}L_p\otimes M_{-p}+c_{p}M_p\otimes L_{-p}+d_{p}M_p\otimes
M_{-p}+e_pY_{p-\frac12}\otimes Y_{\frac12-p}),
\end{eqnarray*}
where all the coefficients are in $\F$ and the sums are all finite.
Since the elements of the form $u_{1,p}:=L_p\otimes
L_{-p}-L_{-p}\otimes L_{p}$, $u_{2,p}:=L_p\otimes
M_{-p}-M_{-p}\otimes L_{p},$ $u_{3,p}:=M_p\otimes
M_{-p}-M_{-p}\otimes M_{p}$ and $u_{4,p}:=Y_{p-\frac12}\otimes
Y_{\frac12-p}-Y_{\frac12-p}\otimes Y_{p-\frac12}$ are all in ${\rm
Im}(1-\tau),$ replacing $v$ by $v-u$, where $u$ is a combination of
some $u_{1,p}$, $u_{2,p}$ ,$u_{3,p}$ and $u_{4,p}$, one can suppose
\begin{eqnarray}\label{wpqr2}
&&c_p=0,\ \forall\ \,p\in\Z;\ \ a_{p},\ d_{p}\ne 0\,\Longrightarrow\
\,p>0\ \mbox{ or }\ p=0;\ \ e_{p}\ne 0\,\Longrightarrow\ \,p>0.
\end{eqnarray}
Then $v$ can be rewritten as
\begin{eqnarray}\label{sm1}
v=\mbox{$\sum\limits_{p\in\Z_{+}}a_{p}$}L_p\otimes
L_{-p}+\mbox{$\sum\limits_{p\in\Z}b_{p}$}L_p\otimes
M_{-p}+\mbox{$\sum\limits_{p\in\Z_{+}}d_{p}$}M_p\otimes
M_{-p}+\mbox{$\sum\limits_{p\in\Z_{>0}}e_{p}$}Y_{p-\frac12}\otimes
Y_{\frac12-p}.
\end{eqnarray}
Assume $a_{p}\ne 0$ for some $p>0$. Choose $q>0$ such that $q\ne p$.
Then $L_{p+q}\otimes L_{-p}$ appears in $L_{q}\cdot v,$ but
(\ref{wpqr2}) implies that the term $L_{-p}\otimes L_{p+q}$ does not
appear in $L_q\cdot v$, a contradiction with the fact that $L_q\cdot
v\in {\rm Im}(1-\tau)$. Then one can suppose $a_p=0,\
\forall\,\,p\in\Z^*$ . Similarly, one can also suppose $d_p=0,\
\forall\,\,p\in\Z^*$ and $e_p=0$, $\forall\,\,p\in\Z$. Then
(\ref{sm1}) becomes
\begin{eqnarray}\label{sm2}
v=\mbox{$\sum\limits_{p\in\Z}$}b_{p}L_p\otimes M_{-p}
+a_{0}L_0\otimes L_{0}+d_{0}M_0\otimes M_{0}\,.
\end{eqnarray}
Recall the fact ${\rm Im}(1-\tau)\subset{\rm Ker}(1+\tau)$ and our
hypothesis $\LL\cdot v\subset{\rm Im}(1-\tau)$, one has
\begin{eqnarray*}
0\!\!\!&=&\!\!\!(1+\tau)L_1\cdot v\\
\!\!\!&=&\!\!\!-2a_{0}(L_1\otimes L_{0}+L_0\otimes
L_{1})\\
&&\!\!\!+\mbox{$\sum\limits_{p\in\Z}$}\big((p-2)b_{p-1}
-pb_{p}\big)L_p\otimes
M_{1-p}+\mbox{$\sum\limits_{p\in\Z}$}\big((p-2)b_{p-1}
-pb_{p}\big)M_{1-p}\otimes L_p\,.
\end{eqnarray*}
Comparing the coefficients, and noting that the set
$\{p\,|\,b_p\ne0\}$ is finite, one gets
\begin{eqnarray*}
&&a_0=b_0+b_1=b_{p}=0,\ \ \forall\,\,p\in\Z\SM\{0,1\}.
\end{eqnarray*}
Then (\ref{sm2}) can be rewritten as
\begin{eqnarray}\label{sm3}
v=b_{1}(L_{1}\otimes M_{-1}-L_0\otimes M_{0})+d_0M_0\otimes M_{0}.
\end{eqnarray}
Observing $(1+\tau)L_2\cdot v=0$, one has $b_1=0$.  Thus the lemma
follows.\QED\vskip5pt

\ni{\it Proof of Theorem \ref{theo} (ii) and (iii)}\ \ \rm\  Let
$(\LL ,[\cdot,\cdot],\D)$ be a Lie bialgebra structure on $\LL$.  By
(\ref{deriv}), (\ref{tr}) and Theorem \ref{theo}(i), $\D=\D_r+D$,
where $r\in \VV\,({\rm mod\,}M_0\otimes M_0)$ and $D\in\DD$.   By
(\ref{cLie-s-s}), ${\rm Im}\,\D\subset{\rm Im}(1-\tau)$, so
$\D_r(L_n)+D(L_n)\in$ Im$(1-\tau)$ for $n\in\Z$, which implies that
$\a+\a^\dag=\gamma+\gamma^\dag=0$. Similarly, $\b+\b^\dag=0$ by the
fact that $\D_r(M_n)+D(M_n)\in\mathrm{Im}(1-\tau)$ for $n\in\Z$.
Thus, $D(\LL)\in \mathrm{Im}(1-\tau)$. So
$\mathrm{Im\,}\D_r\in\mathrm{Im}(1-\tau)$. It follows immediately
from  Lemma \ref{lemma3ll} that $r\in{\rm Im}(1-\tau)\,({\rm
mod\,}M_0\otimes M_0)$, proving the first statement of Theorem
\ref{theo}(ii). If $D\in\DD_1$, one can easily verify that
$(1+\xi+\xi^2)\cdot(1\otimes D)\cdot D=0$ by acting it on generators
of $\LL$, which shows $(\LL,[\cdot,\cdot],D)$ is a Lie bialgebra,
and the proof of  Theorem \ref{theo}(ii) is completed. Theorem
\ref{theo}(iii) follows immediately from (\ref{cLie-s-s}),
Definition \ref{def2} and Lemma \ref{some}.
 \QED

\end{document}